\documentclass[12pt,reqno,a4wide]{amsart}

\usepackage{tikz} 
\usepackage{calrsfs}
\usepackage{mathrsfs}

\usepackage{array}
\usepackage{hyperref}

\usetikzlibrary{decorations.pathreplacing}
\usetikzlibrary{fit}

\usepackage{pgfplots}

\allowdisplaybreaks

\catcode`,\active

\catcode`\,12

\begin{document}
	\bibliographystyle{plain}
	
	%%%%%%%%%%%%%%%%%%%%%%%%%%%%%%%%%%%%%%%%%%%%%%%%%%%%%%%%
	%
	%		Title and paper information
	%
	%%%%%%%%%%%%%%%%%%%%%%%%%%%%%%%%%%%%%%%%%%%%%%%%%%%%%%%%
	
	\title[The generating function of A348410 in OEIS using the diagonal method]
	{The generating function of A348410 in OEIS using the diagonal method and another sequence(A001008) from OEIS}

	\author[H. Prodinger ]{Helmut Prodinger }
	\address{Department of Mathematics, University of Stellenbosch 7602, Stellenbosch, South Africa
		and
		NITheCS (National Institute for
		Theoretical and Computational Sciences), South Africa.}
	\email{warrenham33@gmail.com}


	\begin{abstract}
	$a_n=[x^n](1-x)^{-n}(1-x^2)^{-n}$	is the sequence  A348410 in the Encyclopedia of Integer Sequences. Using a method from
	Hautus and Klarner from 1971 and the software \textsf{Gfun} we find an algebraic equation for the generating function $g(z)=\sum_na_nz^n$.
	
	The second identity uses the \emph{generalized binomial series}, popularized in the textbook \emph{Concrete mathematics}.
	\end{abstract}.

	\maketitle

\section{Introduction}

Let 
\begin{equation*}
a_n=[x^n](1-x)^{-n}(1-x^2)^{-n},
\end{equation*}
which is one form of the sequence A348410 in \cite{OEIS}. The interest of the recent note \cite{tongniu2} was to study the
generating function
\begin{equation*}
g(z)=\sum_{n\ge0}a_nz^n.
\end{equation*}
We offer an approach, originally from Hautus and Klarner \cite{HK} that deserves to be better known.

As a bonus, we discuss another recent identity from the Encyclopedia of Integer Sequences \cite{OEIS}.

\section{Using the diagonalization method}

The first step is to define the bivariate generating function
\begin{equation*}
F(x,y)=\sum_{n\ge0}y^n(1-x)^{-n}(1-x^2)^{-n}=\frac{(1+x)(1-x)^2}{1-x-x^2+x^3-y}=\sum_{m,n\ge0}f_{m,n}x^my^n.
\end{equation*}
Then the function of interest $g(z)$ is given by the \emph{diagonal}:
\begin{equation*}
g(z)=\sum_{n\ge0}f_{n,n}z^n.
\end{equation*}
Hautus and Klarner \cite{HK} invented a method to pull out this diagonal using \emph{residues}.

To find the generating function $g(z)$ of the diagonal, one should consider
\begin{equation*}
F\Bigl(zt,\frac 1t\Bigr)\frac 1t=\frac{(1+tz)(1-tz)^2}{-1+t-z^2t^3-zt^2+z^3t^4}
=\frac{(1+tz)(1-tz)^2}{z^3(t-s_1)(t-s_2)(t-s_3)(t-s_4)}.
\end{equation*}
and further
\begin{equation*}
	\frac1{2\pi i}\oint\frac{(1+tz)(1-tz)^2}{z^3(t-s_1)(t-s_2)(t-s_3)(t-s_4)}dt
\end{equation*}
where the integral is taken over a simple contour containing all the singularities $s(z)$ of
the series such that $s(z)$ stays bounded as $z\to0$.
The function of interest is then the sum of the residues  but only those that lie within a small circle around $t=0$.
Now (we use \textsf{Gfun} \cite{gfun} for the following expansions)
\begin{align*}
s_1&=-\frac 1z-\frac14+\frac 18z-\frac{27}{256}z^2+\dots,\\
s_2&=\frac 1z+\frac12\frac{\sqrt{2}}{\sqrt{z}}-\frac 38+\frac {37}{128}\sqrt{2z}+\dots,\\
s_3&=\frac 1z-\frac12\frac{\sqrt{2}}{\sqrt{z}}-\frac 38-\frac {37}{128}\sqrt{2z}+\dots,\\
s_4&=1+z+3z^2+9z^3+32z^4+119z^5+466z^6+1881z^7+\dots.
\end{align*}
So only $s_4$ qualifies and we end up with
\begin{align*}
g(z)&=\text{Res}\Bigl(F\Bigl(zt,\frac 1t\Bigr)\frac 1t;t=s_4\Bigr)\\&=
\text{Res}\Bigl( \frac{(1+tz)(1-tz)^2}{z^3(t-s_1)(t-s_2)(t-s_3)(t-s_4)};t=s_4\Bigr)\\
&=[(t-s_4)^{-1}]\frac{(1+tz)(1-tz)^2}{z^3(t-s_1)(t-s_2)(t-s_3)(t-s_4)}\\
&=\frac{(1+tz)(1-tz)^2}{z^3(t-s_1)(t-s_2)(t-s_3)}\Big|_{t=s_4}.
\end{align*}
We simplify the denominator by differentiation:
\begin{align*}
z^3(t-s_1)(t-s_2)(t-s_3)\Big|_{t=s_4}&=\frac {\partial}{\partial t}z^3(t-s_1)(t-s_2)(t-s_3)(t-s_4)\Big|_{t=s_4}\\
	&=(1-zt)(1-zt-4z^2t^2)\Big|_{t=s_4}.
\end{align*}
The final answer is thus
\begin{align*}
g(z)&=\frac{(1+s_4z)(1-s_4z)}{1-zs_4-4z^2s_4^2}=\frac{1-s_4^2z^2}{1-zs_4-4z^2s_4^2}\\&=1+z+5 z^2+19 z^3+85 z^4+376 z^5+1715z^6+7890z^7\dots\,.
\end{align*}
Since $s_4$ satisfies an algebraic equation (of order 4), one may rewrite this as
\begin{equation*}
g(z)=-{\frac {40{z}^{3}s_4^{3}+16z^4s_4^{3}+28{z}^{3}s_4^{2}-32z^2s_4^{2}-40zs_4-67{z}^{2}s_4+		32-35z-64{z}^{2}}{256{z}^{2}+107z-32}}.
\end{equation*}
This algebraic function satisfies itself an algebraic equation, which we obtain thanks to  the procedure \textsf{algfuntoalgeq} of \textsf{Gfun}. It is
\begin{equation*}
(256z^2+107z-32)g^4+(-256z^2-107z+32)g^3+12z(8z+3)g^2-4z(4z+1)g+z^2=0.
\end{equation*}
\textsf{Gfun} can also produce a holonomic recursion for the coefficients $a_m$ of $g(z)$:
\begin{align*}
2680&(4m+1)(2m+1)(4m+3)a_m\\&+(2909610+5084335m+2915850m^2+552965m^3)a_{m+1}\\&
+(-2005464-927852m+180504m^2+90732m^3)a_{m+2}\\&-64(m+3)(1759m^2+8778m+12080)a_{m+3}\\&+7168(m+4)(m+3)(2m+7)a_{m+4}=0
\end{align*}

\section{Sequence A001008} 

As an epilogue, we offer our approach towards the identity in the title, which appeared recently as \cite{tongniu3}
\begin{equation*}
H_n=\frac1m\sum_{k=1}^n\frac{(-1)^{k+1}}{k}\binom{mk}{k}\binom{n+(m-1)k}{n-k},
\end{equation*}
with harmonic numbers $H_n$. We start with the RHS and form the generating function
\begin{align*}
\mathcal{S}&=\sum_{n\ge 1}z^n\frac1m\sum_{k=1}^n\frac{(-1)^{k+1}}{k}\binom{mk}{k}\binom{n+(m-1)k}{n-k}\\
&=\frac1m\sum_{k\ge 1}\sum_{j\ge0}z^{k+j}\frac{(-1)^{k+1}}{k}\binom{mk}{k}\binom{j+mk}{j}\\
&=\frac1m\sum_{k\ge 1}z^{k}\frac{(-1)^{k+1}}{k}\binom{mk}{k}\frac{1}{(1-z)^{mk+1}}\\
&=\frac{-1}{m(1-z)}\sum_{k\ge 1}\frac{1}{k}\binom{mk}{k}u^k,\quad \text{with}\quad u:=\frac{-z}{(1-z)^m}.
\end{align*}
Note the $u(1-z)^m=-z$, or with $\mathcal{B}=1-z$, it leads to
\begin{equation*}
u\mathcal{B}^m=\mathcal{B}-1 ,\quad \text{or}\quad \mathcal{B}=1+u\mathcal{B}^m.
\end{equation*}
This function $\mathcal{B}$ is the \emph{generalized binomial series}, as discussed in \cite{GKP} around page 200.
The substitution $-z=w$ is not doing much but makes life a bit easier. Consider
\begin{equation*}
\sum_{k\ge 1}\frac{1}{k}\binom{mk}{k}u^k ,\quad \text{with}\quad u=\frac{w}{(1+w)^m}.
\end{equation*}
We claim that 
\begin{equation*}
	\frac1m\sum_{k\ge 1}\frac{1}{k}\binom{mk}{k}u^k =\log(1+w).
\end{equation*}
In order to get rid of the logarithm we study a differentiated form
\begin{equation*}
	\frac1m\sum_{k\ge 1}\frac{1}{k}\binom{mk}{k}ku^{k-1} =\frac{(1+w)^{m+1}}{1+w-wm}\frac{1}{1+w}.
\end{equation*}
Multiplying it by $u$ leads to an equivalent form
\begin{equation*}
	\frac1m\sum_{k\ge 1}\binom{mk}{k}u^{k} =\frac{w}{1+w-wm}.
\end{equation*}
Addying the term for $k=0$ and then multiplying by $m$,
\begin{equation*}
	\sum_{k\ge 0}\binom{mk}{k}u^{k} =\frac{1+w}{1+w-wm}.
\end{equation*}
The key to the success is the formula (5.61) in \cite{GKN}, which evaluates the sum:
\begin{equation*}
\sum_{k\ge 0}\binom{mk}{k}u^{k}=\frac1{1-m+m\mathcal{B}^{-1}}=\frac1{1-m+\frac{m}{1+w}}=\frac{1+w}{1+w-m-mw+m}	=\frac{1+w}{1+w-wm},
\end{equation*}
as desired.  
Then
\begin{equation*}
\mathcal{S}=-\frac1{1-z}\log(1+w)=-\frac1{1-z}\log(1-z)=\sum_{n\ge1}H_nz^n,
\end{equation*}
finishing the proof.

\end{document}